\documentclass[12pt]{article}
\usepackage{amsmath}
\usepackage[usenames]{color}
\usepackage{mathrsfs}
\usepackage{amsfonts}
\usepackage{amssymb,amsmath}
\usepackage{CJK}
\usepackage{cite}
\usepackage{cases}
\usepackage{amsthm}

\def\d{\delta}
\def\cl{\centerline}

\def\al{\alpha}
\def\vs{\vspace*}

\def\Z{\mathbb{Z}}
\def\C{\mathbb{C}}
\def\F{\mathbb{C}}
\def\QED{\hfill$\Box$}

\numberwithin{equation}{section}
\newtheorem{theo}{Theorem}[section]

\newtheorem{lemm}[theo]{Lemma}
\newtheorem{prop}[theo]{Proposition}
\def\ptl{\partial}
\def\G{\Gamma}
\begin{document}
\begin{CJK*}{GBK}{song}

\begin{center}
{\bf\large Structures of not-finitely graded Lie algebras\\ related to generalized Heisenberg-Virasoro algebras\,$^*$}
\footnote {$^*\,$Supported by the National Natural Science Foundation of China (No. 11371278, 11271284) and the grant no.~12XD1405000 of Shanghai Municipal Science and Technology Commission.

$^{\,\dag}$ Corresponding author: X. Yue (xiaoqingyue@tongji.edu.cn).
}
\end{center}

\cl{Guangzhe Fan$^{\,*}$,  Chenhong Zhou$^{\,*}$, Xiaoqing Yue$^{\,*,\dag}$}

\cl{\small $^{\,*}$ Department of Mathematics, Tongji
University, Shanghai 200092, China}
\cl{\small E-mail:  yzfanguangzhe@126.com,\ zhouchenhong66@163.com,\ xiaoqingyue@tongji.edu.cn}\vs{5pt}

\vs{8pt}

{\small
\parskip .005 truein
\baselineskip 3pt \lineskip 3pt

\noindent{{\bf Abstract:}  In this paper, we study the structure theory of a class of not-finitely graded Lie algebras related to generalized Heisenberg-Virasoro  algebras. In particular, the derivation algebras, the automorphism groups and the second cohomology groups
 of these Lie algebras are determined. \vs{5pt}

}

}
\parskip .001 truein\baselineskip 6pt \lineskip 6pt
\section{Introduction}
\setcounter{equation}{0}

The Heisenberg-Virasoro algebra contains the classical Heisenberg algebra and Virasoro algebra as subalgebras. As the universal central extension of the Lie algebra of differential operators on a circle of order at most one, the Heisenberg-Virasoro algebra is an important object in mathematics and  physics, and has been widely studied. For example, the twisted Heisenberg-Virasoro algebra was firstly studied by E. Arbarello et al. in \cite{A}. Various generalizations of the Heisenberg-Virasoro algebra have been extensively studied by several authors (e.g., \cite{LWZ,LZ,B,SJ}). However, it seems to us that little has been known on not-finitely graded aspect of generalized Heisenberg-Virasoro algebras.

In recent years, some researches on Lie algebras concerning their derivation algebras, automorphisms, second cohomology groups have been undertaken by many authors (see, e.g., \cite{CHS,DZ,F,S1,S2,SXZ,SZ1,SZ2,SZhou,WWY}). It is well known that central extensions, which are determined by second cohomology groups, are closely related to the structures of Lie algebras (see, e.g., \cite{K,KP}). The computation of the second cohomology groups seems to be important and interesting. Not-finitely graded Lie algebras are important objects in Lie theory, whose structure and representation theories are subjects of studies with more challenges than those of finitely graded Lie algebras. In \cite{CHS}, the authors have studied the structure theory of a class of not-finitely graded Lie algebras related to generalized Virasoro algebras.

In this paper, we consider the following Lie algebras, which are referred to as {\it  generalized Heisenberg-Virasoro algebra $HV(\Gamma)$}:
Let $\Gamma$ be any nontrivial additive subgroup of $\C$, and $\C[\Gamma\times\Z_+]$  the semigroup algebra of $\Gamma\times\Z_+$
with basis $\{x^{\al,i}:=x^\al t^i\,|\,\alpha\in \Gamma,\,i\in\Z_+\}$ and product $x^{\al,i}x^{\beta,j}=x^{\al+\beta,i+j}$. Let $\ptl_x,\,\ptl_t$ be the derivations of $\F[\Gamma\times\Z_+]$ defined by
$\ptl_x(x^{\al,i})=\al x^{\al,i}$, $\ptl_t(x^{\al,i})=i x^{\al,i-1}$ for $\al\in\Gamma,\,i\in\Z_+$. Denote $\ptl=\ptl_x+\ptl_t$. Then $HV(\Gamma)$ is the Lie algebra with the underlying space $\F[\Gamma\times\Z_+]\ptl\oplus\F[\G\times\Z_+]$, basis $\{L_{\al,i}:=x^{\al,i}\ptl$, $H_{\beta,j}:=x^{\beta,j}\,|\,\al,\beta\in\G,\,i,j\in\Z_+\}$ and relations
\begin{eqnarray}\label{algebra-rela1}
&\!\!\!\!\!\!\!\!\!\!\!\!\!\!\!\!\!\!&
[L_{\al,i},L_{\beta,j}]=
(\beta-\al)L_{\al+\beta,i+j}+(j-i)L_{\al+\beta,i+j-1}, \\
\label{algebra-rela2}
&\!\!\!\!\!\!\!\!\!\!\!\!\!\!\!\!\!\!&
[L_{\al,i},H_{\beta,j}]=
\beta H_{\al+\beta,i+j}+j H_{\al+\beta,i+j-1},\\
\label{algebra-rela3}&\!\!\!\!\!\!\!\!\!\!\!\!\!\!\!\!\!\!&
[H_{\al,i},H_{\beta,j}]=0.
\end{eqnarray}
Let $W$ be $\F[\G\times\Z_+]\ptl$ with basis
$\{L_{\al,i}:=x^{\al,i}\ptl\,|\,\al\in\G,\,i\in\Z_+\}$ and  relation
\eqref{algebra-rela1},
and  $H$ be $\F[\G\times\Z_+]$ with basis
$\{H_{\beta,j}:=x^{\beta,j}\,|\,\beta\in\G,\,j\in\Z_+\}$ and  relation
\eqref{algebra-rela3}.
We simply denote $HV=HV(\Gamma)$, then $W$ and $H$ are subalgebras of $HV$. The algebra $W$ is the generalized Witt algebra ${\cal W}=W(0,1,0;\G)$ of Witt type studied in \cite{SXZ}.

Let us recall some basic concepts. A Lie algebra $\cal L$ is  {\it finitely graded} if
there exists an abelian group $G$ such that ${\cal L}=\oplus_{a\in G}{\cal L}_{[a]}$ is $G$-graded satisfying
\begin{equation}\label{Fi-g}\
\mbox{$[{\cal L}_{[a]},{\cal L}_{[b]}]\subset {\cal L}_{[a+b]}$ \ and \ ${\rm dim\,}{\cal L}_{[a]}<\infty$ \ for $a,b\in G$.}\end{equation}
An element $x\in {\cal L}$ is {\it ${\rm ad}$-locally finite} if for any $y\in {\cal L}$ the subspace ${\rm span}\{{\rm ad}_x^i(y)\,|\,i\in\Z_+\}$ is finite-dimensional, where ${\rm ad}_x:y\mapsto[x,y]\ (y\in \cal L)$ is the {\it adjoint operator} of $x$.

Similar to the Lie algebras studied in \cite{CHS}, the Lie algebra $HV$  has the following significant features:\begin{itemize}\parskip-3pt\item[\rm(1)]
It  has a {\it finitely graded filtration} in the sense that there exists a filtration~$0\subset{HV}^{(0)}\subset{HV}^{(1)}\subset\cdots,$ satisfying
$[{HV}^{(i)},{HV}^{(j)}]\subseteq {HV}^{(i+j)}$ for all $i,j\in\Z_+$ and
each ${HV}^{(i)}$ is  finitely graded, i.e., there exists some abelian group $G$ which is independent of $i$ (one can simply choose $G=\G$ in this case) such that
${HV}^{(i)}=\oplus_{g\in G}{HV}^{(i)}_g$,
$[{HV}^{(i)}_g,{HV}^{(j)}_h]\subseteq {HV}^{(i+j)}_{g+h},$
${\rm dim\,}{HV}^{(i)}_g<\infty$ for all $g,h\in G,\,i,j\in\Z_+$ (cf.~\eqref{Fi-g}) and ${HV}_0^{(0)}\neq 0$.
\item[\rm(2)]It has the set $\mathfrak{F}$ of ad-locally finite elements, where \begin{eqnarray}
\label{Locally}\mathfrak{F}=\{aL_{0,0}+\mbox{$\sum\limits_{\alpha,j}$}b_{\alpha,j}H_{\alpha,j}\mbox{ (finite~sum)}\,|\, a,b_{\alpha,j}\in\mathbb{C}, ~ \alpha\in\Gamma,j\in \Z_+\}.
\end{eqnarray}

\item[\rm(3)] $HV$ has exactly two proper ideals: $\mathbb{C}H_{0,0}$ and $H$. In particular, $H$ is the unique maximal ideal of $HV$ (cf.~Theorem \ref{theo1}).
\end{itemize}

 In this paper we shall mainly study the structure theory of $HV$ (namely, derivations, automorphisms, 2-cocycles).
 The Lie algebra $HV$ is $\G$-graded
\begin{equation*}\label{a1-0}
HV=\raisebox{-5pt}{${}^{\, \, \displaystyle\oplus}_{\alpha\in \Gamma}$}{HV}_{\alpha},\ \ \
 {HV}_{\alpha}={\rm span} \{L_{\alpha,i}, H_{\alpha,i}\mid i\in \Z_+\}\mbox{ \ for }\al\in\Gamma.
 \end{equation*}
However, it is not finitely graded. Nevertheless, as stated in \cite{CHS}, due to the fact that $\Gamma$ may not be finitely generated (as a group), $HV$ may not be finitely generated as a Lie algebra. The classical techniques (such as those in \cite{F}) cannot be directly applied to our situation here. One must employ
some new techniques in order to tackle  problems associated with not-finitely graded and not-finitely generated Lie algebras (this is also one of our motivations to present our results here). For instance, one of our strategies used in the present paper is to study the actions of  $L_{\alpha,0}, L_{0,i}$ respectively so that the determination of derivation algebra and automorphism group can be done much more efficiently. The main results of the present paper are summarized in Theorems \ref{theo1}, \ref{theo3}, \ref{theo4}, \ref{theo5} and \ref{theo5.1}.

Throughout the paper, we denote by $\C,\,\C^*,\, \Z,\, \Z_+,\,\G^*$ the sets of complex numbers, nonzero complex numbers, integers, nonnegative integers, nonzero elements of $\G$ respectively.

\section{Some properties of $HV$}
Now we study some properties of the Lie algebra $HV$, which will be summarized in Theorem \ref{theo1}. Firstly, we recall the following concept.

Suppose $A$ is a nonzero ideal of $HV$, and
$0\ne x=\sum_{\al\in\G}x_\al\in A$ (finite sum) such that each homogeneous component $x_\al=\sum_{i\in\Z_+}(a_{\al,i}L_{\al,i}+b_{\al,i}H_{\al,i})$ (finite sum) for some $a_{\al,i},b_{\al,i}\in\C$.
Let \begin{equation}\label{GGG-x}{\rm Supp\,}x=\{\al\in\G\,|\,x_{\al}\ne0\},\end{equation}
then it is called the {\it support} of $x$.
\begin{theo}\label{theo1}\begin{itemize}\parskip-3pt\item[\rm(1)]The Lie algebra $HV$ is not finitely graded.
\item[\rm(2)]
HV has the set $\mathfrak{F}$ of ad-locally finite elements, where \begin{eqnarray}
\label{Local}\mathfrak{F}=\{aL_{0,0}+\mbox{$\sum\limits_{\alpha,j}$}b_{\alpha,j}H_{\alpha,j}\mbox{ (finite~sum)}\,|\, a,b_{\alpha,j}\in\mathbb{C}, ~ \alpha\in\Gamma,j\in \Z_+\}.
\end{eqnarray}
\item[\rm(3)] $HV$ has exactly two proper ideals: $\mathbb{C}H_{0,0}$ and $H$. In particular, $H$ is the unique maximal ideal of $HV$.
\end{itemize}
\end{theo}
\noindent{\it Proof~}~(1)~Firstly, we assume that there exists some abelian group $G$ such that $HV=\oplus_{g\in G}HV_{[g]}$ is $G$-graded satisfying \eqref{Fi-g}. Suppose $HV_{[0]}\neq 0$, then we can take any nonzero $x\in HV_{[0]}$ and $\beta\in \G$ such that $\beta\notin {\rm Supp}x$ (cf.~\eqref{GGG-x}). If we assume that $L_{\beta,2}\in \oplus_{g\in I}HV_{[g]}$ for some finite subset $I$ of $G$, we obtain that $\{{\rm ad}^k_x L_{\beta,2}\,|\,k\in\Z_+\}$ is an independent subset of $\oplus_{g\in I}HV_{[g]}$. In particular, $\sum_{g\in I}{\rm dim}HV_{[g]}=\infty$. Since $I$ is finite, there exists some $g_0\in I$ such that dim $HV_{[g_0]}=\infty$, contradicting our assumption that  $HV_{[g]}$ is finite dimensional for any $g\in G$.

So we have $HV_{[0]}=0$. Choose a total order ``$\prec$'' on $G$ compatible with its group structure. Suppose $L_{0,0}=\sum_{i=1}^nx_{g_i}$ with $g_1\prec\cdots \prec g_{n-1}\prec g_n$ for some $n\geq 1$, where $0\neq x_{g_i}\in HV_{[g_i]}$.  Due to
\begin{equation}\label{eq-alpha-------}[L_{0,0}, L_{\alpha, 0}]=\alpha L_{\alpha,0},
\end{equation}$L_{0,0}$ preserves every one-dimensional space $\C L_{\alpha,0}$ ($\alpha\in\G$), we have that $L_{0,0}$ can not lie in $\oplus_{0\prec g\in G}HV_{[g]}$ or $\oplus_{G\ni g\prec 0}HV_{[g]}$ of $HV$. So we have $0\prec g_n$ and $g_1\prec 0$.

By comparing the maximal homogenous components (with respect to the $G$-gradation of $HV$) of both sides of \eqref{eq-alpha-------},
we obtain that the maximal homogenous component of $L_{\alpha,0}$ lies in $HV_{[g_n]}$ for all $\alpha\in\Gamma$. So there exists $\lambda_{\alpha}\in\C^*$ such that $L_{\alpha,0}=\lambda_{\alpha}L_{0,0}+\sum_{g\prec g_n}y^{\alpha}_g$, where $y^\alpha_g\in HV_{[g]}$. In fact, substituting $L_{\alpha,0}=\lambda_{\alpha}L_{0,0}+\sum_{g\prec g_n}y^{\alpha}_g$ on the left hand side of \eqref{eq-alpha-------} and then proceeding the analysis above gives $L_{\alpha,0}=\lambda_{\alpha}L_{0,0}+\sum_{g\prec0 }y^{\alpha}_g$. Similarly, we can deduce $L_{\alpha,0}=\lambda_{\alpha}L_{0,0}+\sum_{0\prec g }y^{\alpha}_g$. This leads to a contradiction.

(2)~It can be obtained easily by the relations of $HV$.

(3)~Let $H'$ be a proper ideal of $HV$ which satisfies $\mathbb{C}H_{0,0}\neq H'$, then there must exist a nonzero element $\sum_{\gamma,s}a_{\gamma,s}L_{\gamma,s}+ \sum_{\beta,j}c_{\beta,j}H_{\beta,j}$ in $H'\setminus\mathbb{C}H_{0,0}$.
 %with $(\gamma,s)\neq(0,0)$.
 Let any $L_{\alpha,i}$ and $H_{\alpha,i}$ act on it respectively, we get
\begin{eqnarray}\label{2.2}\!\!\!\!\!\!\!
[L_{\al,i},\mbox{$\sum\limits_{\gamma,s}$}a_{\gamma,s}L_{\gamma,s}+ \mbox{$\sum\limits_{\beta,j}$}c_{\beta,j}H_{\beta,j}]&\!\!\!=\!\!\!&~\mbox{$\sum\limits_{\gamma,s}$}a_{\gamma,s}((\gamma-\alpha)L_{\alpha+\gamma,i+s}+(s-i)L_{\alpha+\gamma,i+s-1})\nonumber\\&\!\!\!\!\!\!&
+\mbox{$\sum\limits_{\beta,j}$}c_{\beta,j}(\beta H_{\alpha+\beta,i+j}
 +jH_{\alpha+\beta,i+j-1})
\end{eqnarray}
 and
  \begin{equation}\label{2.3}\
[H_{\al,i},\mbox{$\sum\limits_{\gamma,s}$}a_{\gamma,s}L_{\gamma,s}+ \mbox{$\sum\limits_{\beta,j}$}c_{\beta,j}H_{\beta,j}]=
-\mbox{$\sum\limits_{\gamma,s}$}a_{\gamma,s}(\alpha H_{\alpha+\gamma,i+s}+iH_{\alpha+\gamma,i+s-1}).
  \end{equation}

\begin{itemize}\parskip-3pt\item[\rm(i)]If all the $a_{\gamma,s}$ are equal to zero, there must exist some $c_{\beta,j}\neq0$ with $(\beta,j)\neq(0,0)$. From \eqref{2.2}, then we have $H=H'$.
\item[\rm(ii)]If there exist some $a_{\gamma,s}\neq0$ %with $(\gamma,s)\neq(0,0)$
, we can get $H'=HV$ from \eqref{2.2} and \eqref{2.3}. It is a contradiction.
\end{itemize}
From the above discussion, $HV$ has only two proper ideals:~$\mathbb{C}H_{0,0}$ and $H$. Due to $\mathbb{C}H_{0,0}\subsetneqq H$, we get that $H$ is the unique maximal ideal of $HV$.\QED

\section{ Derivation algebra of $HV$}
Recall that a linear map $D:HV\rightarrow HV$ is a {\it derivation} of $HV$
if $D\big([x,y]\big)=[D(x),y]+[x, D(y)]$ for any $x,y\in HV$.
For any  $z\in HV$,
 the adjoint operator ${\rm ad}_z\!:HV\rightarrow HV$ is a derivation, called an {\it inner
 derivation}. Denote by
 ${\rm Der\,}{HV}$
 and ${\rm ad\,}{HV}$ the vector spaces of all derivations and inner derivations respectively. Then the first cohomology group
${H}^{1}(HV,\ HV)\cong{\rm Der\,}{HV}/{\rm ad\,}{HV}$. We say that a derivation $D\in {\rm Der\,}HV$ is of degree $\gamma$, if it satisfies that  $D(HV_{\al})\subset HV_{\al+\gamma}$. Denote by $({\rm Der\,}HV)_{\gamma}$ the space of all derivations of degree $\gamma$.

Denote by ${\rm Hom}_\Z(\Gamma,\C)$ the space of group homomorphisms from $\Gamma$ to $\C$. For each $\phi\in {\rm Hom}_\Z(\Gamma,\C)$, the scalar multiplication $\phi$ by $c\in\C$ is defined by $(c\phi)(\gamma)=c\phi(\gamma)$, thus ${\rm Hom}_\Z(\Gamma,\C)$ is a vector space. Then for each $\phi\in{\rm Hom}_\Z(\Gamma,\C)$, we can define a derivation $D_{\phi}$ as follows:
\begin{equation}\label{d5}
D_{\phi}( L_{\al,i})=\phi(\alpha)L_{\al,i},\ \ D_{\phi}( H_{\al,i})=\phi(\alpha)H_{\al,i}
\ \ \ \mbox{for}\ \ \al\in\Gamma,\ i\in\Z_+.
\end{equation}
We still use ${\rm Hom}_\Z(\Gamma,\F)$ to denote the corresponding subspace of ${\rm Der\,} HV$. In particular, since $\phi_0:\al\mapsto\al$ is in ${\rm Hom}_\Z(\Gamma,\C)$, we have the derivation
\begin{equation}\label{D0====}D_0=D_{\phi_0}:L_{\al,i}\mapsto \al L_{\al,i},\ \ H_{\al,i}\mapsto \al H_{\al,i} \mbox{ \ for \ }\al\in\G,\ i\in\Z_+.\end{equation}

Now we define two derivations of degree $0$, which are obviously not inner derivations of $HV$.
\begin{equation*}\label{OUTER1}
\left\{\begin{array}{llll}D_1(L_{\alpha,i})=\alpha H_{\alpha,i}+i H_{\alpha,i-1}, D_1(H_{\alpha,i})=0,\\[4pt]D_2(L_{\alpha,i})=0,
~~~~~~~~~~~~~~~~~~~ D_2(H_{\alpha,i})=H_{\alpha,i} &\mbox{for \ }\alpha\in\Gamma,\,i\in\Z_+.
\end{array}\right.\end{equation*}

\begin{theo}\label{theo3}
The derivation space of $HV$ can be written as
\begin{equation}\label{THEOM3}
{\rm Der\,}HV =\raisebox{-5pt}{${}^{\, \, \displaystyle\oplus}_{\gamma\in \Gamma}$}({{\rm Der\,}HV})_{\gamma}={\rm ad}HV\oplus{\rm Hom}_\Z(\Gamma,\F)\oplus\mathbb{C}D_1\oplus\mathbb{C}D_2, \end{equation}
where $({\rm Der\,}HV)_\gamma\subset{\rm ad}HV$, if $\gamma\in\Gamma^*$, and $({\rm Der\,}HV)_0=({\rm ad}HV)_0\oplus{\rm Hom}_\Z(\Gamma,\F)\oplus\mathbb{C}D_1\oplus\mathbb{C}D_2$.
\end{theo}

\noindent{\it Proof~}~ It can be obtained by the following lemmas~3.2$\thicksim$3.6 immediately.\QED

\begin{lemm}\label{lemm1}
Every derivation $D\in HV$ can be written as
 \begin{equation}\label{L-1}
D=\mbox{$\sum\limits_{\gamma\in\G}$} D_\gamma,\ \ D_{\gamma}\in({\rm Der\,}HV)_\gamma \end{equation}
 such that only finitely many $D_\gamma(x)\neq0$ for every $x\in HV$ $($such a sum in \eqref{L-1} is called summable$)$.
\end{lemm}
\noindent{\it Proof~}~For a derivation $D\in{\rm Der\,}HV$ and an element $x_\alpha\in HV_\alpha$, assume that $$D(x_\alpha)=\mbox{$\sum\limits
_{\beta\in\Gamma}$} y_\beta.$$ We define $D_\gamma(x_\alpha)=y_{\alpha+\gamma}$. Then a direct computation shows that $D_\gamma$ is a derivation.\QED

\begin{lemm}\label{lemm2}
%Regarding $\ol D\in{\rm Der\,}\L$ as in ${\rm Der\,}\ol\L$, and b
For any $D\in {\rm Der\,}HV$, replacing $D$ with $D-{\rm ad}_u$ for some $u\in HV$, we can suppose \begin{equation}\label{L-007-D}
D(L_{0,0})=0.\end{equation}
\end{lemm}

\noindent{\it Proof~}~Suppose $D(L_{0,0})\!=\!\sum_{\al,j}(a_{\al,j}L_{\al,j}+c_{\al,j}H_{\al,j})\in HV$ for some $a_{\al,j},\ c_{\al,j}\in\C$.
For any $\al\in\G$, we define $b_{\al,j},\ d_{\al,j} \in\C$ inductively on $j\geq0$ by
\begin{equation}\label{definedn11}
b_{\alpha,j}=\left\{\begin{array}{llll}j^{-1}(-a_{\al,j-1}-\alpha b_{\al,j-1})&\mbox{if \ }j\geq1,\\[4pt]0&\mbox{if \ }j=0.
\end{array}\right.\end{equation}
and
\begin{equation}\label{definedn2}
d_{\alpha,j}=\left\{\begin{array}{llll}j^{-1}(-c_{\al,j-1}-\alpha d_{\al,j-1})&\mbox{if \ }j\geq1,\\[4pt]0&\mbox{if \ }j=0.
\end{array}\right.\end{equation}

Take $u=\sum_{\al,j}(b_{\al,j}L_{\al,j}+d_{\al,j}H_{\al,j})$.
Note that $u\in{HV}$, by \eqref{definedn11} and \eqref{definedn2},  \vspace*{-7pt}it gives $$\begin{array}{ll}
  D(L_{0,0})-{\rm ad}_u(L_{0,0})\!\!\!&=\mbox{$\sum\limits_{\al,j}$}(a_{\al,j}L_{\al,j}+c_{\al,j}H_{\al,j})-\mbox{$\sum\limits_{\al,j}$}\
 (-\alpha b_{\al,j}-(j+1)b_{\al,j+1})L_{\al,j}\\[12pt]&\,\,\,\,\,\, -\mbox{$\sum\limits_{\al,j}$}\
 (-\alpha d_{\al,j}-(j+1)d_{\al,j+1})H_{\al,j}
 \\[12pt]&
=0. \end{array}$$\QED

\begin{lemm}\label{lemm3}
If $\gamma\in\Gamma^*$, $D\in({\rm Der\,}HV)_\gamma$ and $D(L_{0,0})=0$, then
\begin{equation}\label{L-00-D}
D=0.\end{equation}
\end{lemm}

\noindent{\it Proof~}~Applying $D$ to $[L_{0,0}, L_{\alpha,0}]=\alpha L_{\alpha,0}$ with $\alpha\in\G^*$, then we have
 \begin{equation}\label{eq3.4}
[L_{0,0}, D(L_{\alpha,0})]=\alpha D(L_{\alpha,0}).\end{equation}
Assume $D(L_{\alpha,0})=\sum_{j\in\Z_+}(b_j L_{\alpha+\gamma,j}+d_j H_{\alpha+\gamma,j}) $, and by \eqref{eq3.4} we have~
\begin{equation}\label{definedn1}
\left\{\begin{array}{llll}\gamma b_j=-(j+1)b_{j+1}  &\mbox{for \ }j\geq0,\\[4pt]\gamma d_j=-(j+1)d_{j+1}&\mbox{for \ }j\geq0.
\end{array}\right.\end{equation}
As~$b_j=0,j\gg0$, $b_j=0$ for any $j\geq0$. Similarly, $d_j=0$ for any $j\geq0$. Thus we have $D(L_{\alpha,0})=0$, for any $\alpha\in\Gamma$.

 Applying $D$ to $[L_{0,0}, L_{0,1}]= L_{0,0}$, then we have
 \begin{equation}\label{eq3.5}
 [L_{0,0}, D(L_{0,1})]=D(L_{0,0}).\end{equation}
 Assume $D(L_{0,j})=\sum_{k\in\Z_+} (f(j,k) L_{\gamma,k}+g(j,k) H_{\gamma,k})$, for $f(j,k)$, $g(j,k) \in{\mathbb{C}}$, and note that $f(0,k)=g(0,k)=0$. By \eqref{eq3.5}, we have
 \begin{equation}\label{definedn1}
\left\{\begin{array}{llll}\gamma f(1,k)=-(k+1)f(1,k+1) &\mbox{for \ }k\geq0,\\[4pt]\gamma g(1,k)=-(k+1)g(1,k+1) &\mbox{for \ }k\geq0.
\end{array}\right.\end{equation}
Then we can get $f(1,k)=g(1,k)=0$ for any $k\in\Z_+$, and thus $D(L_{0,1})=0$.

Now applying $D$ to $[L_{0,0}, H_{\alpha,0}]=\alpha H_{\alpha,0}$, then we have
 \begin{equation}\label{eq3.11}
[L_{0,0}, D(H_{\alpha,0})]=\alpha D(H_{\alpha,0}).\end{equation}
Assume $D(H_{\alpha,0})=\sum_{j\in\Z_+}(e_j L_{\alpha+\gamma,j}+f_j H_{\alpha+\gamma,j}) $, and from \eqref{eq3.11} we obtain~
\begin{equation}\label{definedn1}
\left\{\begin{array}{llll}\gamma e_j=-(j+1)e_{j+1}  &\mbox{for \ }j\geq0,\\[4pt]\gamma f_j=-(j+1)f_{j+1}&\mbox{for \ }j\geq0.
\end{array}\right.\end{equation}
As~$e_j=0,j\gg0$, $e_j=0$ for any $j\geq0$. Similarly, $f_j=0$ for any $j\geq0$. Thus we have $D(H_{\alpha,0})=0$, for any $\alpha\in\Gamma$.

Since $HV$ is generated by $\{L_{\alpha,0},L_{0,1},H_{\alpha,0}\mid\alpha\in\Gamma\}$ from relation \eqref{algebra-rela1}, we get $D=0$ where $D\in({\rm Der\,}HV)_\gamma$, for $\gamma\in\Gamma^*$.\QED

\begin{lemm}\label{lemm4}
Assume $D\in({\rm Der\,}HV)_0$~  and~ $D(L_{0,0})=0$, we have $:$
\begin{itemize}\parskip-3pt
\item[$(i)$]
%For $D$ satisfying the above conditions, we have
$D(L_{\alpha,0})=b_\alpha L_{\alpha,0}+d_\alpha H_{\alpha,0}$, where $b_\alpha\in{\rm Hom}_\Z(\Gamma,\F)$, $d_\alpha=\alpha d_1 $, $d_1\in{\rm \mathbb{C}},\alpha\in\Gamma $.
    %\item[(ii)]For the $D$ satisfying the above conditions, we have $D(L_{0,j})=jf(1,0)L_{0,j-1}+jg(1,0)H_{0,j-1}$ for any$~j\in{\rm Z_{+}} $, with $f(1,0)\in{\rm \mathbb{C}}$, $g(1,0)\in{\rm \mathbb{C}}$.
    \item[$(ii)$]
    %For $D$ satisfying the above conditions, we have
    $D(H_{\alpha,0})=f_\alpha H_{\alpha,0}$, where $f_\alpha=b_\alpha+f_0$, $f_0\in{\rm \mathbb{C}} $, $\alpha\in\Gamma $.
    %\item[(iv)]For the $D$ satisfying the above conditions, we have $D(H_{0,j})=jf(1,0)H_{0,j-1}+g'(0,0)H_{0,j}$ for any$~\alpha\in\Gamma $,with $f(1,0)\in{\rm \mathbb{C}}$, $g'(0,0)\in{\rm \mathbb{C}}$.
\end{itemize}
\end{lemm}

 \noindent{\it Proof~}~For $(i)$, assume $D(L_{\alpha,0})=\sum_{\alpha,j}(b_{\alpha,j} L_{\alpha,j}+d_{\alpha,j} H_{\alpha,j})$, and from \eqref{eq3.4} we have~
\begin{equation}\label{definedn1}
\left\{\begin{array}{llll}b_{\alpha,j}=0  &\mbox{for \ }j\geq1,\\[4pt]d_{\alpha,j}=0 &\mbox{for \ }j\geq1.
\end{array}\right.\end{equation}
Simply denote $b_{\alpha,0}$, $d_{\alpha,0}$ by $b_\alpha$, $d_\alpha$ respectively, we have $D(L_{\alpha,0})=b_\alpha L_{\alpha,0}+d_\alpha H_{\alpha,0}$, for any $\alpha\in\Gamma$.
Applying $D$ to $[L_{\alpha,0}, L_{\beta,0}]=(\beta-\alpha)L_{\alpha+\beta,0}$, we get that
\begin{equation}\label{definedn1}
\left\{\begin{array}{llll}b_\alpha+b_\beta=b_{\alpha+\beta}  &\mbox{for \ }\alpha\neq\beta,\\[4pt]\beta d_\beta-\alpha d_\alpha=(\beta-\alpha)d_{\alpha+\beta} &\mbox{for \ }\alpha,\ \beta\in\Gamma.
\end{array}\right.\end{equation}
As $D(L_{0,0})=b_0 L_{0,0}=0$, we have $b_0=0$. According to \eqref{definedn1}, we also have $-b_{\alpha}=b_{-\alpha}$.
Set $\alpha=\beta$ in \eqref{definedn1}, we obtain that $b_{2\alpha}=b_{(\alpha+\eta)+(\alpha-\eta)}=b_{\alpha+\eta}+b_{\alpha-\eta}=b_{\alpha}+b_{\eta}+b_{\alpha}+b_{-\eta}=2b_{\alpha}$, for any $\eta\in\Gamma\backslash\{\pm\alpha,0\}$. Thus \eqref{definedn1} holds for all $\alpha,\beta\in\Gamma$, which shows that the map $\Phi:\alpha\longmapsto b_{\alpha}$ is an element in ${\rm Hom}_\Z(\Gamma,\F)$. Hence we get that $b_\alpha\in{\rm Hom}_\Z(\Gamma,\F)$.

Suppose that $\Gamma'=c\Gamma$ for some $c\in\mathbb{C}^*$, we have $HV(\Gamma)\cong HV(\Gamma')$. Thus without loss of generality, we can always suppose $1\in\Gamma$. Then we can obtain that $d_\alpha=(d_1-d_0)\alpha+d_0$, for $d_1\in\mathbb{C}$. Due to $d_0=0$, we have $d_\alpha=\alpha d_1$, for $d_1\in\mathbb{C}$.

 For $(ii)$,  assume $D(H_{0,0})=\sum_{j\in\Z_+}(e_{j} L_{0,j}+f_{j} H_{0,j}) $, and apply $D$ to $[L_{0,0}, H_{0,0}]= 0$. Then we have $e_{j}=f_{j}=0$, for any$~j\geq1 $, and we get $D(H_{0,0})=e_{0} L_{0,0}+f_{0} H_{0,0} $. Applying $D$ to $[L_{0,1}, H_{0,0}]= 0$, we have $e_{0}=0$. Thus $D(H_{0,0})=f_{0} H_{0,0}$.
 Assume $D(H_{\alpha,0})=\sum_{\alpha,j}(e_{\alpha,j} L_{\alpha,j}+f_{\alpha,j} H_{\alpha,j}) $, and by \eqref{eq3.11}, we have~
  \begin{equation}\label{definedn7}
\left\{\begin{array}{llll} e_{\alpha,j}=0  &\mbox{for \ }j\geq1,\\[4pt] f_{\alpha,j}=0 &\mbox{for \ }j\geq1.
\end{array}\right.\end{equation}
So $D(H_{\alpha,0})=e_{\alpha,0} L_{\alpha,0}+f_{\alpha,0} H_{\alpha,0}:= e_{\alpha} L_{\alpha,0}+f_{\alpha} H_{\alpha,0}$.  %we also have $D(H_{\alpha,0})=0$, for all $\alpha\in\Gamma$.
Applying $D$ to $[L_{-\alpha,0}, H_{\alpha,0}]= \alpha H_{0,0}$, we have $e_{\alpha}=0$. Thus we can assume $D(H_{\alpha,0})=f_{\alpha} H_{\alpha,0}$. Now applying $D$ to $[L_{\alpha,0}, H_{\beta,0}]= \beta H_{\alpha+\beta,0}$, we have $b_{\alpha}+f_{\beta}=f_{\alpha+\beta}$. Therefore $b_{\alpha}+f_{0}=f_{\alpha}$, with $f_0\in{\rm \mathbb{C}} $ for any $\alpha\in\Gamma$.\QED

If we replace $D'$ by $D-D_\phi$(cf.~\eqref{d5})(note that this replacement does not affect \eqref{L-007-D}), we can suppose $b_\alpha=0$ for all $\alpha\in\Gamma$. Then we have $D'(L_{\alpha,0})=\alpha d_1 H_{\alpha,0}$ and $D'(H_{\alpha,0})=f_0 H_{\alpha,0}$.

Assume $D'(L_{\alpha,1})=\sum_{j\in\Z_+}(f(1,j) L_{\alpha,j}+g(1,j) H_{\alpha,j}) $, for $f(1,j)$, $g(1,j) \in{\mathbb{C}}$. Applying $D'$ to $[L_{0,0}, L_{\alpha,1}]= \alpha L_{\alpha,1}+L_{\alpha,0}$, we get $f(1,j)=g(1,j)=0$ for $j\geq1$. So $D'(L_{\alpha,1})=f(1,0)L_{\alpha,0}+g(1,0)H_{\alpha,0}+\alpha d_1 H_{\alpha,1}$. By induction on $i$, we have $D'(L_{\alpha,i})=if(1,0)L_{\alpha,i-1}+ig(1,0)H_{\alpha,i-1}+\alpha d_1 H_{\alpha,i}$.

Applying $D'$ to $[L_{-\alpha,i}, H_{\alpha,0}]= \alpha H_{0,i}$ for $\alpha\in\G^*$, we have $D'(H_{0,i})=if(1,0)H_{0,i-1}+f_0 H_{0,i}$.
Applying $D'$ to $[L_{0,i}, H_{\alpha,0}]= \alpha H_{\alpha,i}$ for $\alpha\in\G^*$, we obtain $D'(H_{\alpha,i})=if(1,0)H_{\alpha,i-1}+f_0 H_{\alpha,i}$ for $\alpha\in\G^*$. So $D'(H_{\alpha,i})=if(1,0)H_{\alpha,i-1}+f_0 H_{\alpha,i}$ for $\al\in\G,\ i\in\Z_+$. And it follows that
 \begin{equation*}\label{lemma11}
\left\{\begin{array}{llll} D'(L_{\alpha,i})=if(1,0)L_{\alpha,i-1}+ig(1,0) H_{\alpha,i-1}+\alpha d_1 H_{\alpha,i}, \\[4pt] D'(H_{\alpha,i})=if(1,0)H_{\alpha,i-1}+f_0 H_{\alpha,i},
\end{array}\right.\end{equation*}
for any $\al\in\G,\ i\in\Z_+$.
Due to $[L_{0,0}, L_{\alpha,i}]= \alpha L_{\alpha,i}+i L_{\alpha,i-1}$ and $[L_{0,0}, H_{\alpha,i}]= \alpha H_{\alpha,i}+i H_{\alpha,i-1}$, if we replace $u$ (cf.~\eqref{L-007-D})by $u+f(1,0)L_{0,0}$, and denote $D'+D_0$ by $D'$ (cf.~\eqref{D0====}),
then
 \begin{equation}\label{lemma12}
\left\{\begin{array}{llll} D'(L_{\alpha,i})=ig(1,0) H_{\alpha,i-1}+\alpha d_1 H_{\alpha,i}, \\[4pt] D'(H_{\alpha,i})=f_0 H_{\alpha,i},
\end{array}\right.\end{equation}
with $g(1,0),\ f_0, d_1\in\mathbb{C}$, for any $\al\in\G,\ i\in\Z_+$.

Applying $D'$ to $[L_{\al,i},L_{\beta,j}]= (\beta-\al)L_{\al+\beta,i+j}+(j-i)L_{\al+\beta,i+j-1}$, we obtain $d_1=g(1,0)$.
So we have
\begin{equation}\label{lemma13}
\left\{\begin{array}{lllll} D'(L_{\alpha,i})=d_1 (\alpha  H_{\alpha,i}+i H_{\alpha,i-1}), \\[4pt] D'(H_{\alpha,i})=f_0 H_{\alpha,i},
\end{array}\right.\end{equation}
with $f_0, d_1\in\mathbb{C}$, for any $\al\in\G,\ i\in\Z_+$.

\begin{lemm}\label{lemm2}
Under the above notations, we have a decomposition for $({\rm Der\,}HV)_0:$
 \begin{equation}\label{L-2}
({\rm Der\,}HV)_0=({\rm ad}HV)_0\oplus{\rm Hom}_\Z(\Gamma,\F)\oplus\mathbb{C}D_1\oplus\mathbb{C}D_2,
\end{equation}
where $D_1(L_{\alpha,i})=\alpha H_{\alpha,i}+i H_{\alpha,i-1},\ D_1(H_{\alpha,i})=0,\ D_2(L_{\alpha,i})=0,\ D_2(H_{\alpha,i})= H_{\alpha,i}$, for any $\al\in\G,\ i\in\Z_+$.
\end{lemm}

\section{Automorphism group of $HV$}
Firstly, we study the automorphism group of $W$.
 Denote by ${\rm
Aut\,}{W}$
the automorphism group of $W$. Let $\chi(\G)$ be the set of characters of $\G$, i.e., the set of group homomorphisms $\tau:\G\to\C^*$. Set $\G^{\C^*}=\{c\in\C^*\,|\,c\G=\G\}$. We define a group structure on $\chi(\G)\times\G^{\C^*}$ by
\begin{equation}\label{group-sss}
(\tau_1,c_1)\cdot(\tau_2,c_2)=(\tau,c_1c_2),\mbox{ \ where \ }\tau:\al\mapsto\tau_1(c_2\al)\tau_2(\al)\mbox{ for }\al\in\G.
\end{equation}It turns out that the group $\chi(\G)\times\G^{\C^*}$ is just the semidirect product $\chi(\G)\rtimes\G^{\C^*}$ under the action given by $(c\tau)(\alpha)=\tau(c\alpha)$ for all $c\in\G^{\C^*}$, $\tau\in\chi(\G),$  $\alpha\in\Gamma$.
 We define a group homomorphism $\phi:(\tau,c)\mapsto\phi_{\tau,c}$ from $\chi(\G)\times\G^{\C^*}$ to ${\rm Aut\,}W$ such that
 $\phi_{\tau,c}$ is the automorphism of $W$ defined by\begin{equation}\label{Aususus}
 \mbox{$\phi_{\tau,c}:L_{\al,i}\mapsto\tau(\al)c^{i-1}L_{c\al,i}$ \ for \ $\al\in\G,\,i\in\Z_+$.}\end{equation}
One can easily verify that $\phi_{\tau,c}$ is indeed an automorphism of $W$.
\begin{theo}\label{theo4}
${\rm Aut\,}W\cong\chi(\G)\rtimes\G^{\C^*}$.
\end{theo}
\noindent{\it Proof~}~This theorem can be obtained by the following three lemmas.\QED

As $L_{0,0}$ is the unique locally finite element of $W$, for any $\sigma\in{\rm Aut\,}{W}$, we can suppose $\sigma(L_{0,0})=c^{-1}L_{0,0}$, for $c\in\C^*$.
\begin{lemm}\label{Auto-lamm}
For any $\alpha\in \G$, we have
$\sigma(L_{\alpha,0})=c^{-1}\tau(\alpha)L_{c\alpha,0}$, where $c\in\C^*,\ \tau\in{\rm Hom}_\Z(\Gamma,\F^{*})$.
\end{lemm}
\noindent{\it Proof~}~Assume that $\sigma(L_{\alpha,0})= \sum_{\beta,j}b_{\beta,j}^\alpha L_{\beta,j}$ (finite sum) for some $b_{\beta,j}^\alpha\in\C$. Applying $\sigma$ to $[L_{0,0}, L_{\alpha,0}]=\alpha L_{\alpha,0}$, we have
$[L_{0,0},\sum_{\beta,j}b_{\beta,j}^\alpha L_{\beta,j}]=c\alpha\sum_{\beta,j}b_{\beta,j}^\alpha L_{\beta,j}$.
By \eqref{algebra-rela1},
we get $$\mbox{$\sum\limits_{\beta,j}$}b_{\beta,j}^\alpha(\beta L_{\beta,j}+j L_{\beta,j-1})= c\alpha \mbox{$\sum\limits_{\beta,j}$}b_{\beta,j}^\alpha L_{\beta,j}.$$
Comparing the coefficients of $L_{\beta,j}$ on both sides of the above formula, we have $\beta b_{\beta,j}^\alpha+(j+1)b_{\beta,j+1}^\alpha=c\alpha b_{\beta,j}^\alpha$. Then it gives that $\beta=c\alpha,\ b_{\beta,j}^\alpha=0,\ j\geq1$.
So $\sigma(L_{\alpha,0})=\varphi(\alpha)L_{c\alpha,0}=c^{-1}\tau(\alpha)L_{c\alpha,0}$, where
$\varphi$ is a function from ${\rm \Gamma}$ to ${\rm \C^{*}}$ such that $\varphi(0)=c^{-1}$. Let $\tau(\alpha)=c\varphi(\alpha)$. We can verify that $\tau\in{\rm Hom}_\Z(\Gamma,\F^{*})$.\QED

\begin{lemm}\label{Annnclai1}
For any $i\in{\rm \Z_{+}} $, we have $\sigma(L_{0,i})=c^{i-1}L_{0,i}$.
\end{lemm}
\noindent{\it Proof~}~The lemma holds for $i=0$. Now we suppose it holds for $i-1$. Applying $\sigma$ to
$[L_{0,0}, L_{0,i}]=i L_{0,i-1}$, we have $[L_{0,0}, \sigma(L_{0,i})]=i c^{i-1}L_{0,i-1}$. Furthermore, we obtain $\sigma(L_{0,i})=c^{i-1}L_{0,i}+a L_{0,0}$, for $a\in\C$. Now applying $\sigma$ to $[L_{0,1}, L_{0,i}]=(i-1) L_{0,i}$ for $i>1$, we get $a=0$. Then it is obvious that $\sigma(L_{0,i})=c^{i-1}L_{0,i}$.\QED

\begin{lemm}\label{Annnclai}
For any $\alpha\in\Gamma,\ i\in{\rm \Z_{+}} $, we have $\sigma(L_{\alpha,i})=\tau(\alpha)c^{i-1}L_{c\alpha,i}$.
\end{lemm}
\noindent{\it Proof~}~If $i=0$, the lemma holds. Now we suppose it is true for $i-1$. Applying $\sigma$ to
$[L_{\alpha,0}, L_{0,i}]=-\alpha L_{\alpha,i}+iL_{\alpha,i-1}$ with $\alpha\in\G^*$ and using $\sigma(L_{\alpha,i-1})=\tau(\alpha)c^{i-2}L_{c\alpha,i-1}$, we can deduce $\sigma(L_{\alpha,i})=\tau(\alpha)c^{i-1}L_{c\alpha,i}$, for any $\alpha\in\Gamma^*,\ i\in{\rm \Z_{+}} $. Then, by Lemma \ref{Annnclai1}, we have $\sigma(L_{\alpha,i})=\tau(\alpha)c^{i-1}L_{c\alpha,i}$, for any $\alpha\in\Gamma,\ i\in{\rm \Z_{+}} $.\QED

Next we compute the automorphism group of $H$.

\begin{lemm}\label{clai4}
For any $\alpha\in \Gamma$, we have
$\sigma(H_{\alpha,0})=\tau(\alpha)eH_{c\alpha,0}$, where $e\in\mathbb{C}^{*},\ \tau\in{\rm Hom}_\Z(\Gamma,\F)$.
\end{lemm}
\noindent{\it Proof~}~Applying $\sigma$ to $[L_{0,0}, H_{\alpha,0}]=\alpha H_{\alpha,0}$, since $H$ is an abelain ideal$(cf.(\ref{algebra-rela2}), (\ref{algebra-rela3}))$, we get $[L_{0,0}, \sigma(H_{\alpha,0})]=c \alpha\sigma(H_{\alpha,0})$. That is to say, $\sigma(H_{\alpha,0})$ is an eigenvector of ${\rm ad}_{L_{0,0}}$ with eigenvector $c\alpha$. Hence we obtain that $\sigma(H_{\alpha,0})=f(\alpha)H_{c\alpha,0}$. Furthermore applying $\sigma$ to $[L_{\alpha,0}, H_{\beta,0}]=\beta H_{\alpha+\beta,0}$ gives $\tau(\alpha)c^{-1}[L_{c\alpha,0}, f(\beta)H_{c\beta,0}]=\beta f(\alpha+\beta) H_{c(\alpha+\beta),0}$, i.e., $f(\beta)\tau(\alpha)=f(\alpha+\beta)$. Denote $f(0)=e$, we get $f(\alpha)=\tau(\alpha)e$, where $\tau\in{\rm Hom}_\Z(\Gamma,\F)$. Due to Theorem 2.1(3), then we have $e\in\mathbb{C}^{*}$.\QED

\begin{lemm}\label{clai5}
For any $j\in {\rm \Z_{+}}$, we have
$\sigma(H_{0,j})=ec^{j}H_{0,j}$.
\end{lemm}
\noindent{\it Proof~}~Applying $\sigma$ to $[L_{0,1}, H_{0,j}]=j H_{0,j}$, and we get $[L_{0,1}, \sigma(H_{0,j})]=j\sigma(H_{0,j})$. Assume that $\sigma(H_{0,j})=\sum_{\beta,k}x_{\beta,k} H_{\beta,k}$, we have $\sum_{\beta,k}x_{\beta,k}(\beta H_{\beta,k+1}+k H_{\beta,k})= j\sum_{\beta,j}x_{\beta,k}H_{\beta,k}$.
Comparing the coefficients of $H_{\beta,k}$ gives $\beta x_{\beta,{k-1}}=(j-k)x_{\beta,k}$. Thus $\beta=0$. We can assume that $\sigma(H_{0,j})=e_{j} H_{0,j},\,e_j\in{\mathbb{C}^*}$. Applying $\sigma$ to $[L_{0,0}, H_{0,j}]=j H_{0,j-1}$, we get $e_j=ce_{j-1},\ j>0$. As $e_0=e$, we have $e_j=e c^j $. Thus the lemma holds.\QED

\begin{lemm}\label{clai6}
For any $\alpha\in\Gamma,i\in{\rm Z_{+}} $, we have $\sigma(H_{\alpha,j})=\tau(\alpha)ec^j H_{\alpha,j}$.
\end{lemm}
\noindent{\it Proof~}~If $\alpha=0$, the lemma holds. If $\alpha\neq0$, by applying $\sigma$ to $[L_{0,j}, H_{\alpha,0}]=\alpha H_{\alpha,j}$, we have $[c^{j-1}L_{0,j},\tau(\alpha)e H_{c\beta,0}]=\alpha c^j e \tau(\alpha)  H_{c\alpha,j}=\alpha \sigma(H_{\alpha,j})$. Thus we get $\sigma(H_{\alpha,j})=\tau(\alpha)ec^j H_{\alpha,j}$.\QED

Now we study the automorphism group of $HV$. Denote by ${\rm
Aut\,}{HV}$
the automorphism group of $HV$. Let $\sigma\in{\rm
Aut\,}{HV}$. We have $\sigma(H)=H,\ \sigma(H_{0,0})=e H_{0,0}$, for $e\in{\rm\mathbb{C^*}}$. Let $\overline{\sigma}\in{\rm
Aut\,}{HV/H}(i.e.,\,{\rm
Aut\,}{W})$, and we have $\overline{\sigma}\overline{(L_{\alpha,i})}=\tau(\alpha)c^{i-1}\overline{(L_{\alpha,i})}$, i.e., $\sigma(L_{\alpha,i})\equiv\tau(\alpha)c^{i-1}L_{c\alpha,i}$~mod $H$.

Since $\sigma(L_{0,0})\equiv c^{-1}L_{0,0}$~mod $H$, we can assume $\sigma(L_{0,0})=c^{-1}L_{0,0}+\sum_{\alpha,i}a_{\alpha,i}H_{\alpha,i}$ for some $a_{\alpha,i}\in\C$. Define
\begin{equation}\label{definedn31}
b_{\alpha,i}=\left\{\begin{array}{llll}-i^{-1}a_{\al,i-1}&\mbox{for \ }i>0,\\[4pt]0&\mbox{for \ }i=0.
\end{array}\right.\end{equation}
Denote $\tau=e^{{\rm ad}\sum_{\alpha,i}b_{\alpha,i}H_{\alpha,i}} \in{\rm Inn}HV $, so $\tau(c^{-1}L_{0,0})=c^{-1}L_{0,0}+\sum_{\alpha,i}b_{\alpha,i}H_{\alpha,i}$. Let $\pi=\tau^{-1}\sigma$. We obtain $\pi(L_{0,0})=\tau^{-1}\sigma(L_{0,0})=c^{-1}L_{0,0}$. Then by Lemma \ref{clai4},\ \ref{clai5} and \ref{clai6}, we have $\pi(L_{\alpha,i})=\tau(\alpha)c^{i-1}L_{c\alpha,i}$.
Thus we complete the proof of the following theorem.
\begin{theo}\label{theo5}
${\rm Aut\,}HV\cong{\rm Inn}HV \rtimes\Big((\chi(\G)\rtimes\G^{\C^*})\rtimes{\mathbb{C}^*}\Big)$.
\end{theo}

\section{ Second cohomology group of $HV$}
Recall that a bilinear form $\psi:HV\times HV
\rightarrow \F$ is called a {\it 2-cocycle} on $HV$ if
the following conditions are satisfied:
\begin{eqnarray*}
&&\psi(x,y)=-\psi(y,x),%\\ &&
\ \ \ \psi(x,[y,z])+\psi(y,[z,x])+\psi(z,[x,y])=0,
\end{eqnarray*}
for $x,y,z\in HV$. Denote by
$C^{2}(HV,\ \F)$ the vector space of 2-cocycles on
 $HV$. For any linear function $f: HV \rightarrow
 \F$, one can define a 2-cocycle $\psi_{f}$ by
$ % \begin{eqnarray*}&&
\psi_{f}(x,y)=f([x,y])$ for $x, y \in HV.$ %\end{eqnarray*}
Such a 2-cocycle is called a {\it trivial 2-cocycle} or a {\it 2-coboundary} on $HV$. Denote
by $B^{2}(HV,\ \F)$ the vector space of 2-coboundaries
on $HV$. The quotient space
 \begin{eqnarray*}
&&H^{2}(HV,\
\F)=C^{2}(HV,\ \F)/B^{2}(HV,\ \F)
\end{eqnarray*}
is called  the {\it 2-cohomology group} of $HV$.
There exists a one-to-one correspondence between the set of
equivalence classes of one-dimensional central extensions of
$HV$ by $\F$ and the 2-cohomology group of
$HV$.

\begin{theo}\label{theo5.1}
%We have ${\rm dim}H^2(HV,\mathbb{C})=3$, where
We have $H^2(HV,\C)=0$.
\end{theo}
\noindent{\it Proof~}~It can be
%proved by the reference.
easily obtained by \cite{SZ1}, Proposition \ref{theo5.2} and Proposition \ref{theo5.3}.\QED

\begin{prop}\label{theo5.2}We have
\begin{equation}\label{2-cococo7}
~~~~\phi(L_{\al,i},H_{\beta,j})=0,
\end{equation}
with $\alpha,\beta\in\Gamma$, $i,j\in\Z_+$.
 \end{prop}
\noindent{\it Proof~}~Without loss of generality, we can always suppose $1\in\G$. Let $\psi\in C^2(HV,\C)$, we define a linear function $f:HV\to\C$ such that $f(L_{\al,i})$ and $f(H_{\al,i})$ are defined by induction on $i$ as follows
\begin{equation}\label{f==}
f(L_{\al,i})=\left\{\begin{array}{lll}
\frac1{i+1}\psi(L_{0,0},L_{0,i+1})&\mbox{if \ }\alpha=0,\\[4pt]
\frac1\alpha\Big(\psi(L_{0,0},L_{\alpha,i})-if(L_{\alpha,i-1})\Big)&\mbox{if \ }\alpha\neq0.
\end{array}\right.
\end{equation}
\begin{equation}\label{f====}
f(H_{\al,i})=\left\{\begin{array}{lll}
\psi(L_{0,0},H_{0,1})&\mbox{if \ }\alpha=i=0,\\[4pt]
\frac1i\psi(L_{0,1},H_{0,i})&\mbox{if \ }\alpha=0,\,i\neq0,\\[4pt]
\frac1\alpha\Big(\psi(L_{0,0},H_{\alpha,i})-if(H_{\alpha,i-1})\Big)&\mbox{if \ }\alpha\neq0.
\end{array}\right.
\end{equation}Set $\phi=\psi-\psi_f$, we have
\begin{eqnarray}
\label{phi===++222+1}
\phi(L_{0,0},H_{\beta,j})&\!\!\!=\!\!\!&\psi(L_{0,0},H_{\beta,j})\!-\!f([L_{0,0},H_{\beta,j}])\!=\!0
\mbox{ \ for \ } \beta\neq0.
\end{eqnarray}
And we can easily verify that $\phi(L_{0,0},H_{0,0})=\phi(L_{0,0},H_{0,1})=0$.~Furthermore, for $j\neq1$ and $j\neq0$, we obtain
\begin{eqnarray*}\!\!\!\!\!\!\!\!\!\!\!\!&\!\!\!\!\!\!\!\!\!\!\!\!\!\!\!&
\ \ \ \ \ \ \ \ \phi(L_{0,0},H_{0,j})=\psi(L_{0,0},H_{0,j})-jf(H_{0,j-1})
\nonumber\\\!\!\!\!\!\!\!\!\!\!\!\!&\!\!\!\!\!\!\!\!\!\!\!\!\!\!\!&
\ \ \ \ \ \ \ \ \ \ \ \ \ \ \ \ \ \ \ \ \ \ \ \ \,=\psi(L_{0,0},H_{0,j})-\frac{j}{j-1}\phi(L_{0,1},H_{0,j-1})
\nonumber\\\!\!\!\!\!\!\!\!\!\!\!\!&\!\!\!\!\!\!\!\!\!\!\!\!\!\!\!&\ \ \ \ \ \ \ \ \ \ \ \ \ \ \ \ \ \ \ \ \ \ \ \ \,
=\psi(L_{0,0},H_{0,j})-\frac{1}{j-1}\phi(L_{0,1},[L_{0,0},H_{0,j}])
\nonumber\\\!\!\!\!\!\!\!\!\!\!\!\!&\!\!\!\!\!\!\!\!\!\!\!\!\!\!\!&\ \ \ \ \ \ \ \ \ \ \ \ \ \ \ \ \ \ \ \ \ \ \ \ \,
=0.
\end{eqnarray*}
%$\phi_2(L_{0,0},H_{0,j})=\psi(L_{0,0},H_{0,j})-jf(H_{0,j-1})=\psi(L_{0,0},H_{0,j})-\frac{j}{j-1}\phi_2(L_{0,1},H_{0,j-1})=\psi(L_{0,0},H_{0,j})-\frac{1}{j-1}\phi_2(L_{0,1},[L_{0,0},H_{0,j}])=0$.
Thus
\begin{equation}
\label{phi===1}
\phi(L_{0,0},H_{\beta,j})=0 \ \ \ \ \ \ \mbox{ for } \ \beta\in\Gamma,\ j\in\Z_+.
\end{equation}

Now we want to prove $\phi(L_{\alpha,i},H_{0,0})=0$, for any $\alpha\in\Gamma, i\in\Z_+$.
If $\alpha\neq0$, we use induction on $i$. It holds for $i=0$. Suppose it also holds for $i-1$, then we get
\begin{eqnarray}
\label{phi=1110}
\phi(L_{\alpha,i},H_{0,0})&\!\!\!=\!\!\!&\frac{1}{\alpha}\psi([L_{0,0},L_{\alpha,i}]\!-\!iL_{\alpha,i-1},H_{0,0})\!=\!0\mbox{ \ for \ } \ \alpha\neq0.
\end{eqnarray}
%So $\phi_2(L_{\alpha,i},H_{0,0})=0$ for $\alpha\neq0$.
We can also have
\begin{eqnarray}
\label{phi=1112}
\phi(L_{0,i},H_{0,0})&\!\!\!=\!\!\!&\psi(\frac{1}{i+1}[L_{0,0},L_{0,i+1}],H_{0,0})\!=\!0
\mbox{ \ for \ } i\neq0.
\end{eqnarray}
Therefore, by $(\ref{phi=1110})$ and $(\ref{phi=1112})$, it follows that \begin{equation}
\label{phi===2}
\phi(L_{\alpha,i},H_{0,0})=0\ \ \ \ \ \mbox{ for } \ \alpha\in\Gamma,\ i\in\Z_+.
\end{equation}

Similarly, for $\alpha+\beta\neq0$, we obtain
\begin{eqnarray*}\!\!\!\!\!\!\!\!\!\!\!\!&\!\!\!\!\!\!\!\!\!\!\!\!\!\!\!&
\ \ \ \ \ \ \ \ \phi(L_{\alpha,0},H_{\beta,0})=\psi(L_{\alpha,0},H_{\beta,0})-\beta f(H_{\alpha+\beta,0})
\nonumber\\\!\!\!\!\!\!\!\!\!\!\!\!&\!\!\!\!\!\!\!\!\!\!\!\!\!\!\!&
\ \ \ \ \ \ \ \ \ \ \ \ \ \ \ \ \ \ \ \ \ \ \ \ \ \,
=\psi(L_{\alpha,0},H_{\beta,0})- \frac{\beta}{\alpha+\beta}\psi(L_{0,0},H_{\alpha+\beta,0})
\nonumber\\\!\!\!\!\!\!\!\!\!\!\!\!&\!\!\!\!\!\!\!\!\!\!\!\!\!\!\!&
\ \ \ \ \ \ \ \ \ \ \ \ \ \ \ \ \ \ \ \ \ \ \ \ \ \,
=\psi(L_{\alpha,0},H_{\beta,0})-\frac{1}{\alpha+\beta}\psi(L_{0,0},[L_{\alpha,0},H_{\beta,0}])
\nonumber\\\!\!\!\!\!\!\!\!\!\!\!\!&\!\!\!\!\!\!\!\!\!\!\!\!\!\!\!&
\ \ \ \ \ \ \ \ \ \ \  \ \ \ \ \ \ \ \ \ \ \ \ \ \ \,
=0.
\end{eqnarray*}
%$\phi_2(L_{\alpha,0},H_{\beta,0})=\psi(L_{\alpha,0},H_{\beta,0})-\beta f(H_{\alpha+\beta,0})=\psi(L_{\alpha,0},H_{\beta,0})- \frac{\beta}{\alpha+\beta}\psi(L_{0,0},H_{\alpha+\beta,0})=\psi(L_{\alpha,0},H_{\beta,0})-\frac{1}{\alpha+\beta}\psi(L_{0,0},[L_{\alpha,0},H_{\beta,0}])=0$, for $\alpha+\beta\neq0$.
By induction on $i+j$, for $\alpha+\beta\neq0$, we have
\begin{eqnarray*}\!\!\!\!\!\!\!\!\!\!\!\!&\!\!\!\!\!\!\!\!\!\!\!\!\!\!\!&
\ \ \ \ \ \ \ \ \phi(L_{\alpha,i},H_{\beta,j})=\psi(L_{\alpha,i},H_{\beta,j})-\beta f(H_{\alpha+\beta,i+j})-jf(H_{\alpha+\beta,i+j-1})
\nonumber\\\!\!\!\!\!\!\!\!\!\!\!\!&\!\!\!\!\!\!\!\!\!\!\!\!\!\!\!&
\ \ \ \ \ \ \ \ \ \ \ \ \ \ \ \ \ \ \ \ \ \ \ \ \,=\psi(L_{\alpha,i},H_{\beta,j})-\frac{\beta}{\alpha+\beta}(\psi(L_{0,0},H_{\alpha+\beta,i+j})-(i+j)f(H_{\alpha+\beta,i+j-1}))
\nonumber\\\!\!\!\!\!\!\!\!\!\!\!\!&\!\!\!\!\!\!\!\!\!\!\!\!\!\!\!&
\ \ \ \ \ \ \ \ \ \ \ \ \ \ \ \ \ \ \ \ \ \ \ \ 
\ \ \ \ -\frac{j}{\alpha+\beta}(\psi(L_{0,0},H_{\alpha+\beta,i+j-1})-(i+j-1)f(H_{\alpha+\beta,i+j-2}))
\nonumber\\\!\!\!\!\!\!\!\!\!\!\!\!&\!\!\!\!\!\!\!\!\!\!\!\!\!\!\!&
\ \ \ \ \ \ \ \ \ \ \ \ \ \ \ \ \ \ \ \ \ \ \ \ \,
=-\frac{1}{\alpha+\beta}(j\phi(L_{\alpha,i},H_{\beta,j-1})
+i\phi(L_{\alpha,i-1},H_{\beta,j}))
\nonumber\\\!\!\!\!\!\!\!\!\!\!\!\!&\!\!\!\!\!\!\!\!\!\!\!\!\!\!\!&
\ \ \ \ \ \ \ \ \ \ \ \ \ \ \ \ \ \ \ \ \ \ \ \ \,
=0.
\end{eqnarray*}
Thus we can conclude that
\begin{equation}
\label{phi===2}\phi(L_{\al,i},H_{\beta,j})=\d_{\al+\beta,0}\phi(L_{\alpha,i},H_{-\alpha,j})\ \ \ \ \ \mbox{ for } \alpha,\beta\in\Gamma,\ i,j\in\Z_+.
\end{equation}

From the above discussion, we have already known $\phi(L_{0,0},H_{0,1})=\phi(L_{0,1},H_{0,0})=0$. Now for the case $i+j>1$, we have
\begin{eqnarray*}\!\!\!\!\!\!\!\!\!\!\!\!&\!\!\!\!\!\!\!\!\!\!\!\!\!\!\!&
\ \ \ \ \ \ \ \ \phi(L_{0,i},H_{0,j})=\psi(L_{0,i},H_{0,j})-j f(H_{0,i+j-1})
\nonumber\\\!\!\!\!\!\!\!\!\!\!\!\!&\!\!\!\!\!\!\!\!\!\!\!\!\!\!\!&
\ \ \ \ \ \ \ \ \ \ \ \ \ \ \ \ \ \ \ \ \ \ \ \ \,=\psi(L_{0,i},H_{0,j})-\frac{j}{i+j-1}\psi(L_{0,1},H_{0,i+j-1})
\nonumber\\\!\!\!\!\!\!\!\!\!\!\!\!&\!\!\!\!\!\!\!\!\!\!\!\!\!\!\!&
\ \ \ \ \ \ \ \ \ \ \ \ \ \ \ \ \ \ \ \ \ \ \ \ \,
=\psi(L_{0,i},H_{0,j})-\frac{1}{i+j-1}\psi(L_{0,1},[L_{0,i},H_{0,j}])
\nonumber\\\!\!\!\!\!\!\!\!\!\!\!\!&\!\!\!\!\!\!\!\!\!\!\!\!\!\!\!&
\ \ \ \ \ \ \ \ \ \ \ \ \ \ \ \ \ \ \ \ \ \ \ \ \,
=0.
\end{eqnarray*}
%$\phi_2(L_{0,i},H_{0,j})=\psi(L_{0,i},H_{0,j})-j f(H_{0,i+j-1})=\psi(L_{0,i},H_{0,j})-\frac{j}{i+j-1}\psi(L_{0,1},H_{0,i+j-1})=\psi(L_{0,i},H_{0,j})-\frac{1}{i+j-1}\psi(L_{0,1},[L_{0,i},H_{0,j}])=0$.
So we obtain
\begin{equation}
\label{phi===3}
\phi(L_{0,i},H_{0,j})=0\ \ \ \ \ \mbox{ for } i,j\in\Z_+.
\end{equation}

In order to complete the proof, by $(\ref{phi===2})$ and $(\ref{phi===3})$, we need to compute $\phi(L_{\alpha,i},H_{-\alpha,j})$ for $\alpha\neq0$. If $i+j=1$, we have $\phi(L_{\alpha,0},H_{-\alpha,1})=\phi(L_{\alpha,1},H_{-\alpha,0})=0$, for $\alpha\neq0$. For $i+j>1$, by induction on $i+j$, it follows that
\begin{eqnarray*}\!\!\!\!\!\!\!\!\!\!\!\!&\!\!\!\!\!\!\!\!\!\!\!\!\!\!\!&
\ \ \ \ \ \ \ \ \phi(L_{\alpha,i},H_{-\alpha,j})=\psi(L_{\alpha,i},H_{-\alpha,j})-\frac{1}{i+j}\psi(L_{0,1},[L_{\alpha,i},H_{-\alpha,j}]-jH_{0,i+j-1})-j f(H_{0,i+j-1})=0.
\end{eqnarray*}
%$\phi_2(L_{\alpha,i},H_{-\alpha,j})=\psi(L_{\alpha,i},H_{-\alpha,j})-\frac{1}{i+j}\psi(L_{0,1},[L_{\alpha,i},H_{-\alpha,j}]-jH_{0,i+j-1})-j f(H_{0,i+j-1})=0$.
Thus
 \begin{equation}
\label{phi===3ss}
\phi(L_{\alpha,i},H_{-\alpha,j})=0\ \ \ \ \ \mbox{ for } \alpha\neq0,\ i+j\geq1.
\end{equation}
Based on the above discussion, we obtain
\begin{equation}
\label{phi===2'}\phi(L_{\al,i},H_{\beta,j})=\d_{\al+\beta,0}\d_{i+j,0}\phi(L_{\alpha,0},H_{-\alpha,0})\ \ \ \ \ \mbox{ for } \ \alpha,\beta\in\Gamma,\ i,j\in\Z_+.
\end{equation}

Finally, we consider $\phi(L_{\alpha,0},H_{-\alpha,0})$. Obviously, we have
\begin{eqnarray*}\!\!\!\!\!\!\!\!\!\!\!\!&\!\!\!\!\!\!\!\!\!\!\!\!\!\!\!&
\ \ \ \ \ \ \ \ \phi(L_{\alpha,0},H_{-\alpha,0})=\psi(L_{\alpha,0},H_{-\alpha,0})- f([L_{\alpha,0},H_{-\alpha,0}])
\nonumber\\\!\!\!\!\!\!\!\!\!\!\!\!&\!\!\!\!\!\!\!\!\!\!\!\!\!\!\!&
\ \ \ \ \ \ \ \ \ \ \ \ \ \ \ \ \ \ \ \ \ \ \ \ \ \ \,
=\psi(L_{\alpha,0},H_{-\alpha,0})+\alpha \psi(L_{0,0},H_{0,1})
\nonumber\\\!\!\!\!\!\!\!\!\!\!\!\!&\!\!\!\!\!\!\!\!\!\!\!\!\!\!\!&
\ \ \ \ \ \ \ \ \ \ \ \ \ \ \ \ \ \ \ \ \ \ \ \ \ \ \,
\triangleq c(\alpha).
\end{eqnarray*}
Since $\phi([L_{\alpha,0},L_{\beta,0}],H_{-\alpha-\beta,0})+\phi([L_{\beta,0},H_{-\alpha-\beta,0}],L_{\alpha,0})+\phi([H_{-\alpha-\beta,0},L_{\alpha,0}],L_{\beta,0})=0$, we get
\begin{equation}
\label{phi===4}
(\beta-\alpha)c(\alpha+\beta)+(\beta+\alpha)c(\alpha)=(\beta+\alpha)c(\beta), \mbox{\ \ for\ any\ \ } \alpha,\beta\in\Gamma.
\end{equation}
In \eqref{phi===4}, replacing $\alpha$ by $\alpha-1$ and taking $\beta=1,2$ respectively, we have
\begin{equation}\label{f====1}
\left\{\begin{array}{lll}
(2-\alpha) c(\alpha)+\alpha c(\alpha-1)=\alpha c(1),\\[4pt]
(3-\alpha) c(\alpha+1)+(\alpha+1)c(\alpha-1)=(\alpha+1)c(2).
\end{array}\right.
\end{equation}
And by $(1-\alpha) c(\alpha+1)+(\alpha+1) c(\alpha)=(\alpha+1) c(1)$, we obtain
 \begin{equation}
\label{phi===5}
c(\alpha)=\frac{c(2)}{2}(\alpha^2-\alpha)\ \ \ \ \ \mbox{ for } \alpha\neq0.
\end{equation}
As $\phi(L_{0,0},H_{0,0})=0=c(0)$ also satisfies \eqref{phi===5}, we can conclude that \eqref{phi===5} holds for any $\alpha\in\Gamma$, where $c(2)\in\mathbb{C}$.
So $\phi(L_{\al,0},H_{-\al,0})=b({\al^2-\al})$, for $b\in\C$.

Finally, using the equality $\phi(L_{0,0},[L_{\alpha,1},H_{-\alpha,0}])=0$, we get $\phi(L_{\al,0},H_{-\alpha,0})=0$.\QED

\begin{prop}\label{theo5.3}We have
\begin{equation}\label{2-cococo}
\phi(H_{\al,i},H_{\beta,j})=0,
\end{equation}
with $\alpha,\beta\in\Gamma$, $i,j\in\Z_+$.
\end{prop}
\noindent{\it Proof~}~Let $\psi\in C^2(HV,\C)$, we define a linear function $f:HV\to\C$ (cf.~\eqref{f==} and \eqref{f====})~and set $\phi=\psi-\psi_f$. Then we obtain
\begin{eqnarray}
\label{phi===++222+1}
\phi(H_{0,0},H_{\beta,0})&\!\!\!=\!\!\!&\psi(H_{0,0},H_{\beta,0})\!-\!f([H_{0,0},H_{\beta,0}])\!=\!0.
\end{eqnarray}
By induction on $j$, we have
\begin{equation}
\label{phi===}
\phi(H_{0,0},H_{\beta,j})=0\mbox{ \ for \ } \beta\in\Gamma,\,j\in\Z_+.
\end{equation}
Furthermore, for $\al\neq0$, we obtain
\begin{eqnarray*}\!\!\!\!\!\!\!\!\!\!\!\!&\!\!\!\!\!\!\!\!\!\!\!\!\!\!\!&
\ \ \ \ \ \ \ \ \phi(H_{\al,0},H_{\beta,0})=\psi(H_{\alpha,0},H_{\beta,0})-f([H_{\alpha,0},H_{\beta,0}])
\nonumber\\\!\!\!\!\!\!\!\!\!\!\!\!&\!\!\!\!\!\!\!\!\!\!\!\!\!\!\!&
\ \ \ \ \ \ \ \ \ \ \ \ \ \ \ \ \ \ \ \ \ \ \ \ \ \,=\frac{1}{\alpha}\psi([L_{0,0},H_{\alpha,0}],H_{\beta,0})
\nonumber\\\!\!\!\!\!\!\!\!\!\!\!\!&\!\!\!\!\!\!\!\!\!\!\!\!\!\!\!&
\ \ \ \ \ \ \ \ \ \ \ \ \ \ \ \ \ \ \ \ \ \ \ \ \ \,=-\frac{\beta}{\alpha}\phi(H_{\al,0},H_{\beta,0}).
\end{eqnarray*}
Thus $\phi(H_{\al,0},H_{\beta,0})=0$ for $\alpha+\beta\neq0$.
%\begin{equation}
%\label{phi===++1}
%\phi_3(H_{\al,0},H_{\beta,0})=\psi(H_{\alpha,0},H_{\beta,0})-f([H_{\alpha,0},H_{\beta,0}])=-\frac{\beta}{\alpha}\phi_3(H_{\al,0},H_{\beta,0})\mbox{ \ for \ }\al\neq0.
%\end{equation}

Now considering the case of $\alpha=0$. For $\beta\neq0$, we also have
\begin{eqnarray*}\!\!\!\!\!\!\!\!\!\!\!\!&\!\!\!\!\!\!\!\!\!\!\!\!\!\!\!&
\ \ \ \ \ \ \ \ \phi(H_{0,i},H_{\beta,0})=\psi(H_{0,i},H_{\beta,0})-f([H_{0,i},H_{\beta,0}])
\nonumber\\\!\!\!\!\!\!\!\!\!\!\!\!&\!\!\!\!\!\!\!\!\!\!\!\!\!\!\!&
\ \ \ \ \ \ \ \ \ \ \ \ \ \ \ \ \ \ \ \ \ \ \ \ \,=-\frac{(-1)^i i!}{\beta^i}\phi(H_{0,0},H_{\beta,0})
\nonumber\\\!\!\!\!\!\!\!\!\!\!\!\!&\!\!\!\!\!\!\!\!\!\!\!\!\!\!\!&
\ \ \ \ \ \ \ \ \ \ \ \ \ \ \ \ \ \ \ \ \ \ \ \ \,
=0.
\end{eqnarray*}
By induction on $j$, we get
\begin{equation}
\label{5.30}\phi(H_{0,i},H_{\beta,j})=0 \mbox{ \ for \ }\ \beta\neq0.
\end{equation}

Next we want to prove $\phi(H_{\al,i},H_{\beta,j})=0$ for $\alpha\neq0$ and $i+j\geq1$. Obviously, it holds for $i=0,\ j=1$ and $i=1,\ j=0$. Now assume $\phi(H_{\al,i-1},H_{\beta,j})=\phi(H_{\al,i},H_{\beta,j-1})=0$.
For $\alpha\neq0$, considering
\begin{eqnarray*}\!\!\!\!\!\!\!\!\!\!\!\!&\!\!\!\!\!\!\!\!\!\!\!\!\!\!\!&
\ \ \ \ \ \ \ \ \phi(H_{\al,i},H_{\beta,j})=\psi(H_{\alpha,i},H_{\beta,j})-f([H_{\alpha,i},H_{\beta,j}])
\nonumber\\\!\!\!\!\!\!\!\!\!\!\!\!&\!\!\!\!\!\!\!\!\!\!\!\!\!\!\!&
\ \ \ \ \ \ \ \ \ \ \ \ \ \ \ \ \ \ \ \ \ \ \ \ \ \,=-\frac{\beta}{\alpha}\phi(H_{\al,i},H_{\beta,j})-\frac{j}{\alpha}\phi(H_{\al,i},H_{\beta,j-1})
-\frac{i}{\alpha}\phi(H_{\al,i-1},H_{\beta,j}),
\end{eqnarray*}
%$\phi_3(H_{\al,i},H_{\beta,j})=\psi(H_{\alpha,i},H_{\beta,j})-f([H_{\alpha,i},H_{\beta,j}])
%=-\frac{\beta}{\alpha}\phi(H_{\al,i},H_{\beta,j})-\frac{j}{\alpha}\phi_3(H_{\al,i},\\H_{\beta,j-1})-\frac{i}{\alpha}\phi(H_{\al,i-1},H_{\beta,j})$ for $\alpha\neq0$, we get $\phi_3(H_{\al,i},H_{\beta,j})=-\frac{\beta}{\alpha}\phi_3(H_{\al,i},H_{\beta,j})$.
it follows that
\begin{equation}
\label{phi===++8}
\phi(H_{\al,i},H_{\beta,j})=0 \mbox{ \ for \ }\alpha+\beta\neq0, \ \ i+j\geq1.
\end{equation}
Therefore, we obtain
\begin{equation}
\label{phi===++9}
\phi(H_{\al,i},H_{\beta,j})=\d_{\al+\beta,0}\phi(H_{\al,i},H_{-\alpha,j})\mbox{ \ for \ }\alpha\neq0.
\end{equation}

It is easy to see that $\phi(H_{0,i},H_{0,j})=0 $, for $i, j\in \Z_+$.
Then we have
\begin{equation*}
\label{phi===++3}
\phi(H_{\al,i},H_{-\alpha,j})=-\frac{\alpha}{i-1}\phi(H_{0,j+1},H_{-\alpha,i})
 +\frac{j}{i-1}\phi(H_{0,j},H_{-\alpha,i})=0\mbox{ \ for \ }\ \alpha\neq0,\ i\geq2.
\end{equation*}
%$\phi_3(H_{\al,i},H_{-\alpha,j})=-\frac{\alpha}{i-1}\phi_3(H_{0,j+1},H_{-\alpha,i})
% +\frac{j}{i-1}\phi_3(H_{0,j},H_{-\alpha,i})=0$,
%for any $\alpha\neq0,\ i\geq2$.
Since $\phi(H_{\al,0},H_{-\alpha,1})=0$ and $\phi(H_{\al,1},H_{-\alpha,0})=0$, we get $\phi(H_{\al,i},H_{-\alpha,j})=0$  for $\alpha\neq0,\ i\geq1$. It follows that $\phi(H_{\al,i},H_{-\alpha,j})=\d_{i+j,0}\phi(H_{\al,0},H_{-\alpha,0})$.
Thus we have
\begin{equation}
\label{phi===++9}
\phi(H_{\al,i},H_{\beta,j})=\d_{\al+\beta,0}\d_{i+j,0}\phi(H_{\al,0},H_{-\alpha,0})\ \ \ \ \ \mbox{ for } \ \alpha,\beta\in\Gamma,\ i,j\in\Z_+.
\end{equation}
Now suppose $\phi(H_{\al,0},H_{-\alpha,0})={g(\alpha)}$ and assume $1\in\Gamma$, we get $\phi(H_{1,0},H_{-1,0})={g(1)}\in\mathbb{C}$
and
%\begin{eqnarray*}\!\!\!\!\!\!\!\!\!\!\!\!&\!\!\!\!\!\!\!\!\!\!\!\!\!\!\!&
%\ \ \ \ \ \ \ \ \phi_3(H_{\al,0},H_{-\alpha,0})=\psi([L_{\al-1,0},H_{1,0}],H_{-\alpha,0})
%\nonumber\\\!\!\!\!\!\!\!\!\!\!\!\!&\!\!\!\!\!\!\!\!\!\!\!\!\!\!\!&
%\ \ \ \ \ \ \ \ \ \ \ \ \ \ \ \ \ \ \ \ \ \ \ \ \ \ \ \ \ \,=\alpha\psi(H_{1,0},H_{-1,0})
%\nonumber\\\!\!\!\!\!\!\!\!\!\!\!\!&\!\!\!\!\!\!\!\!\!\!\!\!\!\!\!&
%\ \ \ \ \ \ \ \ \ \ \ \ \ \ \ \ \ \ \ \ \ \ \ \ \ \ \ \ \ \,
%=\alpha g(1).
%\end{eqnarray*}
$\phi(H_{\al,0},H_{-\alpha,0})=\psi([L_{\al-1,0},H_{1,0}],H_{-\alpha,0})=\alpha\psi(H_{1,0},H_{-1,0})=\alpha g(1)$, for $g(1)\in\mathbb{C}$.

Finally, using the equality $\phi(L_{0,0},[H_{\alpha,1},H_{-\alpha,0}])=0$, we get $\phi(H_{\al,0},H_{-\alpha,0})=0$.\QED

%\centerline{\includegraphics[scale=1.2]{actmark.eps}}
%\centerline{\small Figure 1\quad Journal mark}

\end{CJK*}
\end{document}